\documentclass[10pt]{amsart}
\usepackage[normalem]{ulem}

\usepackage{filecontents}

\usepackage[sorting=ynt, backend=biber, datamodel=mrnumber, style=numeric-comp, isbn=false, sortcites, maxbibnames=99, giveninits=true]{biblatex}

\assignrefcontextentries*[sorting=ynt]{*}

\usepackage{amsmath,amssymb,amsthm,geometry,bm}
\usepackage[normalem]{ulem}
\usepackage[pdftex]{graphicx}
\usepackage[colorlinks=true, pdfstartview=FitH, linkcolor=blue,
citecolor=blue, urlcolor=blue]{hyperref}
\newfont{\footsc}{cmcsc10 at 8truept}
\newfont{\footbf}{cmbx10 at 8truept}
\newfont{\footrm}{cmr10 at 10truept}

\newcommand{\C}{{\mathbb C}}
\newcommand{\N}{{\mathbb N}}
\newcommand{\Q}{{\mathbb Q}}
\newcommand{\R}{{\mathbb R}}
\newcommand{\Z}{{\mathbb Z}}
\newcommand{\ep}{\varepsilon}

\newcommand{\cA}{{\mathcal A}}

\newcommand{\cN}{{\mathcal N}}
\newcommand{\cR}{{\mathcal R}}
\newcommand{\cS}{{\mathcal S}}

\newcommand{\fA}{{\mathfrak A}}

\renewcommand{\mod}[1]{{\ifmmode\text{\rm\ (mod~$#1$)}\else\discretionary{}{}{\hbox{ }}\rm(mod~$#1$)\fi}}
\newcommand{\li}{\mathop{\rm li}}
\newcommand{\Li}{\mathop{\rm Li}}

\newcommand{\cM}{{\mathcal M}}

\newcommand{\Deltadot}{\mathord{\mathring\Delta}}
\newcommand{\Edot}{\mathord{\mathring E}}

\numberwithin{equation}{section}

\usepackage{verbatim}

\title[Annotated bibliography for comparative prime number theory]{An annotated bibliography for \\ comparative prime number theory}
\author[Martin, Yang, Bahrini, Bajpai, Benl\.{i}, Downey, Li, Liang, Parvardi, Simpson, White, and Yip]{Greg Martin, Pu Justin Scarfy Yang, Aram Bahrini, Prajeet Bajpai, \\ K\"ubra Benl\.{i}, Jenna Downey, Yuan Yuan Li, Xiaoxuan Liang,\\ Amir Parvardi, Reginald Simpson, Ethan Patrick White, and Chi Hoi Yip}
\address{University of British Columbia\\ Department of Mathematics \\ Room 121\\ 1984 Mathematics Road\\Vancouver, BC Canada V6T 1Z2}
\email{gerg@math.ubc.ca}
\subjclass[2010]{11N13 (11Y35)}

\begin{document}
\nocite{*}

\begin{abstract}
The goal of this annotated bibliography is to record every publication on the topic of comparative prime number theory together with a summary of its results. We use a unified system of notation for the quantities being studied and for the hypotheses under which results are obtained.
\end{abstract}
\maketitle

\large
\section{Introduction}\label{intro section}

Comparative prime number theory is the study of number-theoretic quantities, such as functions that count primes with particular properties, and how they compare to one another. It certainly includes (but is not limited to) ``prime number races'', which examine inequalities between the counting functions of primes in arithmetic progressions to the same modulus; indeed, Chebyshev observing the apparent preponderance of primes of the form $4k+3$ over those of the form $4k+1$ was the historical beginning of comparative prime number theory. Studying inequalities between two functions can be rephrased as studying the sign of their difference, and so the methods of comparative prime number theory also extend to studying the sign (and changes of sign) of other number-theoretic quantities that are less directly related to prime-counting functions.

The phrase ``comparative prime number theory'' goes back at least as far as the title of a long sequence of papers of Knapowski and Tur\'an, beginning with~\cite{MR0146156}. That paper begins with a list of several questions that can be interpreted as an attempt to define the scope of the field, as does the first paper~\cite{MR0162771} in a sequel series by the same authors. Other surveys of these topics include papers by Kaczorowski~\cite{MR1645811} and by Ford and Konyagin~\cite{MR1985941}, as well as an expository introduction to the field by Granville and the first author~\cite{MR2202918}.

This being said, there is no ironclad definition of what is and is not comparative prime number theory. Most quantities in this field have ``explicit formulas'' that express them as sums of oscillatory functions indexed by the zeros of $L$-functions of some type (including the Riemann zeta-function). As such, suitably normalized versions of these quantities are expected to have limiting (logarithmic) distribution functions, which are measures that record the frequencies with which the normalized quantities take values in various intervals in the limit (``continuous histograms'' of their values). In our view, the existence of such a limiting distribution is one of the main criteria for deciding whether a topic does or does not belong to the field of comparative prime number theory.

The purpose of this annotated bibliography is to provide a single exhaustive resource that lists every publication in the field of comparative prime number theory, and provides a summary of the results of each publication included. Like any human endeavour, the fulfillment of that goal will be imperfect. More specifically, we have aimed for completeness for all publications through 2023,
as well as an incomplete list of sources from 2024.

The publications in comparative prime number theory over the 170 years of its existence have understandably used a wide variety of notations for the same objects. Another purpose of this work is to propose a unified system of notation for referring to the functions and quantities that are the main objects of study in comparative prime number theory, as well as uniform terminology for the assumptions on zeros of $L$-functions that arise repeatedly when trying to prove theorems about these quantities. In particular, in our summaries of each publication, we have translated the results into this modern unified notation whenever possible, rather than preserving the notation used by the authors. In this respect, this work is more of a scientific resource than a historical document, although of course we hope it has some utility in the latter role (and we have included authors' exact words on a few occasions, particularly when problems or conjectures were first proposed).

Section~\ref{notation section} is therefore a long section presenting this system of notation for elementary functions, prime counting functions and other summatory functions of number-theoretic quantities, their error terms (both normalized and unnormalized), weighted and averaged versions of these quantities, analogues of these quantities over number fields and function fields, functions that count the number of sign changes of these quantities, and (natural and logarithmic) limiting densities and limiting distribution functions. Section~\ref{other notation section}, beginning on page~\pageref{other notation section}, describes objects and theorems that frequently arise in this field, such as Dirichlet characters and $L$-functions, Landau's theorem, explicit formulas, the power-sum method, $k$-functions, and various hypotheses on the zeros of $L$-functions.
Section~\ref{questions section}, beginning on page~\pageref{questions section}, enumerates the types of questions that comparative number theory studies about the quantities from Section~\ref{notation section}.
The annotated bibliography proper begins on page~\pageref{first bib}.

The origin of this manuscript was a literature survey project by the first two authors in 2012; since then, the other authors have contributed significantly and have greatly expanded the extent of this bibliography and the accompanying material.



\section{Notation related to number theory and real analysis}\label{notation section}

We use $\N$ to denote the set of positive integers, and similarly $\Z$, $\R$, and~$\C$ to denote the sets of integers, real numbers, and complex numbers, respectively. We reserve the letter $p$ to denote prime numbers, and sums and products such as $\sum_p$ and $\prod_{p\mid q}$ are restricted to prime values of~$p$.

We use the following standard conventions regarding magnitudes of complex-valued functions~$f$ and~$g$, real-valued functions~$h$, and nonnegative real-valued functions~$r$ and~$s$ (of a complex or real argument~$z$):
\begin{itemize}
    \item $f(z) \ll s(z)$ (due to Vinogradov) means that there exists a constant $C>0$ such that $|f(z)| \le Cs(z)$ for all values of~$z$ under consideration;
    \item $O(s(z))$ (due to Bachmann) represents an unspecified function $f(z)$ with the property that $f(z)\ll s(z)$;
    \item $r(z) \asymp s(z)$ (due to Hardy) means that both $r(z) \ll s(z)$ and $s(z) \ll r(z)$ are true;
    \item $f(z) \sim g(z)$ (also due to Hardy) means that $\lim f(z)/g(z)=1$, where the location of the limit is taken from context (often as $z\to\infty$ through real numbers);
    \item $f(z) = o(s(z))$ (due to Landau) means that $\lim f(z)/s(z)=0$;
    \item $f(z) = \Omega(s(z))$ (due to Hardy and Littlewood) is the negation of $f(z) = o(s(z))$, or equivalently the statement $\limsup |f(z)|/s(z)>0$;
    \item $h(z) = \Omega_+(s(z))$ and $h(z) = \Omega_-(s(z))$ (due in this form to Landau) mean, respectively, that $\limsup h(z)/s(z)>0$ and $\liminf h(z)/s(z)<0$, either of which implies $h(z) = \Omega(s(z))$;
    \item $h(z) = \Omega_\pm(s(z))$ means that both $h(z) = \Omega_+(s(z))$ and $h(z) = \Omega_-(s(z))$ are true.
\end{itemize}

\subsection{Elementary functions}

As is standard in number theory, we use $\phi(n)$ to denote the Euler totient function, which is the number of reduced residue classes modulo~$n$. We use $\omega(n)$ to denote the number of distinct prime factors of~$n$ and $\Omega(n)$ to denote the number of prime factors of $n$ counted with multiplicity. We let $\mu(n)$ and $\Lambda(n)$ denote the M\"obius and von Mangoldt functions, respectively:
\[
\mu(n) = \begin{cases}
(-1)^{\omega(n)}, &\text{if $n$ is squarefree}, \\ 0, &\text{otherwise;}
\end{cases}
\qquad
\Lambda(n) = \begin{cases}
\log p, &\text{if $n=p^r$ for some $r\in\N$}, \\ 0, &\text{otherwise.} 
\end{cases}
\]

We use $(a,q)$ as a shorthand for $\gcd(a,q)$.
For any $(a,q)=1$, we define
\[
c_q(a) = \#\{b\mod q\colon b^2\equiv a\mod q\}
\]
to be the number of ``square roots'' of $a$ modulo~$q$. For brevity we write $c_q=c_q(1)$, which is also the number of real Dirichlet characters\mod q, or equivalently the index $\bigl[ (\Z/q\Z)^\times : \bigl( (\Z/q\Z)^\times \bigr)^2 \bigr]$; it turns out that $c_q=2^{\omega(q)+\eta}$ where $\eta\in\{-1,0,1\}$ depends upon the power of $2$ dividing~$n$. For $(a,q)=1$, it is the case that $c_q(a)$ equals $c_q$ if $a$ is a square\mod q and $0$ otherwise. (Many sources define $c(q,a)$ to be $c_q(a)-1$, which is more convenient for some purposes and less convenient for others.)

We define two closely related logarithmic integrals
\begin{align*}
\li(x) &= \lim_{\ep\to0+} \bigg( \int_0^{1-\ep} \frac{dt}{\log t} + \int_{1+\ep}^x \frac{dt}{\log t} \bigg) = \sum_{k=1}^K \frac{(k-1)!x}{(\log x)^k} + O_K\bigg( \frac x{(\log x)^{K+1}} \bigg) \\
\Li(x) &= \int_2^x \frac{dt}{\log t} = \li(x) - \li(2) \approx \li(x) - 1.04516378.
\end{align*}

\subsection{Prime counting functions}

We use the standard notation for the prime counting functions
\[
\begin{split}
\pi(x) &= \# \{ p \le x \} = \sum_{p\le x} 1 \\
\Pi(x) &= \sum_{n\le x} \frac{\Lambda(n)}{\log n} = \sum_{p^k \le x} \frac1k = \sum_{k=1}^\infty \frac{\pi(x^{1/k})}k \\
\theta(x) &= \sum_{p \le x} \log p \\
\psi(x) &= \sum_{n\le x} \Lambda(n) = \sum_{p^k \le x} \log p = \sum_{p\le x} \bigg\lfloor \frac{\log x}{\log p} \bigg\rfloor \log p = \sum_{k=1}^\infty \frac{\theta(x^{1/k})}k.
\end{split}
\]
We may replace the cutoff variable $x$ with any set $S$ of real numbers, so that for example
\[
\psi(S) = \sum_{n\in S} \Lambda(n)
\quad\text{and}\quad
\Pi\big( (0,x] \big) = \Pi(x)
\quad\text{and}\quad
\theta\big( (x,y] \big) = \theta(y)-\theta(x).
\]

All of these functions have analogues for prime powers restricted to arithmetic progressions:
\begin{align*}
\pi(x;q,a) &= \# \{ p \le x \colon p\equiv a\mod q \} = \sum_{\substack{p\le x \\ p\equiv a\mod q}} 1 \\
\Pi(x;q,a) &= \sum_{\substack{n\le x \\ n\equiv a\mod q}} \frac{\Lambda(n)}{\log n} = \sum_{\substack{p^k\le x \\ p^k\equiv a\mod q}} \frac1k \\
\theta(x;q,a) &= \sum_{\substack{p\le x \\ p\equiv a\mod q}} \log p \\
\psi(x;q,a) &= \sum_{\substack{n\le x \\ n\equiv a\mod q}} \Lambda(n) = \sum_{\substack{p^k\le x \\ p^k\equiv a\mod q}} \log p.
\end{align*}
These counting functions are interesting only in the case $(a,q)=1$, a restriction that we will usually not state explicitly.
Here too we may replace the first argument with a set, so that for example $\pi(S;q,a) = \#\{p\in S\colon p\equiv a\mod q\}$.

When the third argument is a set rather than an integer, the function counts prime powers that are congruent modulo~$q$ to any element of that set; for example, $\theta(x;q,\{1,2\}) = \theta(x;q,1) + \theta(x;q,2)$. In this context,~$\cR$ and~$\cN$ always refer to the quadratic residues and nonresidues, respectively, among the reduced residues modulo~$q$, so that for example
\begin{align*}
\pi(x;q,\cR) &= \# \{ p \le x \colon p \text{ is a quadratic residue} \mod q \} \\
\pi(x;q,\cN) &= \# \{ p \le x \colon p \text{ is a quadratic nonresidue} \mod q \}.
\end{align*}
Note that $\cR$ contains $\phi(q)/c_q$ residue classes\mod q and $\cN$ contains the other $\phi(q)(1-1/{c_q})$ residue classes. We use $\cA = \cN\cup\cR$ to refer to the set of all reduced residue classes.

When these prime counting functions for arithmetic progressions appear with four arguments instead of three, the function is the difference of the counts for the two indicated arithmetic progressions; for example, $\psi(x;q,a,1) = \psi(x;q,a) - \psi(x;q,1)$.
(Warning: some authors use $\Delta$ for differences of this type, but we give a different meaning to $\Delta$ below in Section~\ref{errors section}.)
When these final two arguments are sets, we make the convention that the two counting functions being subtracted are individually normalized by the number of distinct reduced residue classes in each set; for example,
\begin{align}
\theta(x;7,\{1,2\},\{3,4,5\}) &= \tfrac12\theta(x;7,\{1,2\}) - \tfrac13\theta(x;7,\{3,4,5\}) \notag \\ \Pi(x;q,\cN,\cR) &= \frac1{\phi(q)-c_q} \Pi(x;q,\cN) - \frac1{c_q} \Pi(x;q,\cR). \label{N minus R}
\end{align}
(This convention is consistent with the four-argument notation when the last two arguments are single integers, although these is some dissonance between this convention and the three-argument notation when the last argument is a set, since that function is not normalized in this way.)
There is no need for the notation to admit the possibility of two different moduli, since such a difference can always be written using residue classes of the least common multiple of the moduli: for example,
\[
\pi(x;8,1) - \pi(x;5,2) = 4\pi\big( x;40, \{1,9,17,33\}, \{7,17,27,37\} \big) = 3\pi\big( x;40, \{1,9,33\}, \{7,27,37\} \big).
\]

The residue class~$1\mod q$ is special in some ways, and it is thus helpful to define the notation
\[
\pi(x;q,1,\max) = \pi(x;q,1) - \max_{\substack{a\in(\Z/q\Z)^\times \\ a\not\equiv1\mod q}} \pi(x;q,a),
\quad
\pi(x;q,1,\min) = \pi(x;q,1) - \min_{\substack{a\in(\Z/q\Z)^\times \\ a\not\equiv1\mod q}} \pi(x;q,a),
\]
and similarly for other prime counting functions.

\subsection{Prime ideal classes} \label{number fields sec}

For any number field~$K$ (finite extension of~$\Q$), we say that $\alpha\in K$ is totally positive if~$\alpha$ maps to a positive real number under all embeddings of~$K$ in~$\C$.
We call ideals~$\mathfrak{a}$ and~$\mathfrak{b}$ of a number field~$K$ congruent modulo another ideal $\mathfrak{f}\subset K$ if both~$\mathfrak{a}$ and~$\mathfrak{b}$ are coprime to $\mathfrak{f}$ and there exist totally positive algebraic integers $\alpha$ and $\beta$ in $K$ with $\alpha\equiv\beta \equiv 1\mod {\mathfrak{f}}$ such that $\alpha\mathfrak{a}=\beta \mathfrak{b}$. The equivalence classes of ideals modulo~$\mathfrak{f}$ form a group under ideal multiplication, with the principal ideal class~$\mathfrak{K}_0$ as its identity element. For a character~$\chi$ of this group, we abuse notation slightly by defining $\chi(\mathfrak{a})$ on ideals~$\mathfrak{a}$ directly: if $\mathfrak{a}$ is coprime to $\mathfrak{f}$ then we set $\chi(\mathfrak{a})=\chi([\mathfrak{a}])$, where $[\mathfrak{a}]$ is the ideal class\mod {\mathfrak{f}} containing $\mathfrak{a}$, and if $\mathfrak{a}$ is not coprime to $\mathfrak{f}$ then we set $\chi(\mathfrak{a})=0$. We can now define the Hecke--Landau zeta-function $\zeta(s,\chi)$ to be the Dirichlet series
\[
\zeta(s,\chi)=\sum_{\mathfrak{a}} \frac{\chi(\mathfrak{a})}{\mathfrak{Na}^s},
\]
Finally, for an ideal class $\mathfrak{K}$, we define prime ideal counting functions such as
\[
\pi(x,\mathfrak{K})= \sum_{\substack{\mathfrak{Np}\leq x\\\mathfrak{p}\in\mathfrak{K}\\\mathfrak{p}\,\textbf{prime ideal}}} 1, \qquad
\psi(x,\mathfrak{K})= \sum_{\substack{\mathfrak{Np}^m\leq x\\\mathfrak{p}^m\in\mathfrak{K}\\\mathfrak{p}\,\textbf{prime ideal}}} \log \mathfrak{Np}.
\]

\subsection{Error terms for prime counting functions} \label{errors section}

These prime counting functions have well-known main terms, and it is useful to have a standard notation to refer to the error terms obtained by subtracting these main terms, as well as normalized versions of such error terms.
We use $\Delta$ to denote error terms for the standard prime counting functions:
\[
\Delta^\psi(x) = \psi(x) - x, \quad \Delta^\theta(x) = \theta(x) - x, \quad \Delta^\Pi(x) = \Pi(x) - \li(x), \quad \Delta^\pi(x) = \pi(x) - \li(x).
\]
(In this document's article summaries, we will use the above normalizations even when an article subtracts a slightly different main term: we do not distinguish here between $\li(x)$ and $\Li(x)$ and $\sum_{2\le n\le x} 1/\log n$, for example.)
We also use $E$ for normalized versions of these error terms:
\[
E^\psi(x) = \frac{\Delta^\psi(x)}{\sqrt x}, \quad E^\theta(x) = \frac{\Delta^\theta(x)}{\sqrt x}, \quad E^\Pi(x) = \frac{\Delta^\Pi(x)}{\sqrt x/\log x}, \quad E^\pi(x) = \frac{\Delta^\pi(x)}{\sqrt x/\log x}.
\]
While there is not a formula for starting with a general function $f$ and determining the correct denominator to use when defining $E^f$, the normalization factor chosen is the one for which the resulting~$E$ function is expected to have a limiting logarithmic distribution.

It's not uncommon to integrate these error terms: for a function~$f$ such as~$\pi$, $\Pi$, $\theta$, or~$\psi$, we define $\fA^f_0(x) = \Delta^f(x)$ and, for $m\ge1$,
\[
\fA^f_m(x) = \int_0^x \fA^f_{m-1}(t)\,dt.
\]
(Again, we ignore the fact that some articles might use a different lower endpoint for such integrals.) This operation has a predictable effect on summatory functions and explicit formulas: for example, $\fA_m^\psi(x) = \sum_{n\le x} (\Lambda(n)-1)(x-n)^m/m!$ has an explicit formula containing terms of the form $x^{\rho+m}/\rho(\rho+1)\cdots(\rho+m)$. For repeated integration of the absolute error, we also define $\fA^f_{|0|}(x) = |\fA^f(x)|$ and, for $m\ge1$,
\[
\fA^f_{|m|}(x) = \int_0^x \fA^f_{|m-1|}(t)\,dt.
\]
There are similar logarithmic integration operators: we define $A^f_0(x) = \Delta^f(x)$ and, for $m\ge1$,
\[
A^f_m(x) = \int_0^x A^f_{m-1}(t)\,\frac{dt}t.
\]
This operation also predictably affects summatory functions and explicit formulas: for example, $A_m^\psi(x) = \sum_{n\le x} (\Lambda(n)-1)(\log\frac xn)^m$ has an explicit formula containing terms of the form $x^{\rho}/\rho^{m+1}$. We also use the notation $A^f_{|m|}(x)$ for repeated logarithmic integration of the absolute error.

When we count primes in arithmetic progressions, the error terms $\Delta$ include a factor of $\phi(q)$ for simplicity: for example,
\[
\Delta^\psi(x;q,a) = \phi(q)\psi(x;q,a) - x \quad\text{and}\quad \Delta^\pi(x;q,a) = \phi(q)\pi(x;q,a) - \li(x).
\]
The normalized error terms $E$ are then derived from these $\Delta$ as before: for example,
\[
E^\psi(x;q,a) = \frac{\Delta^\psi(x;q,a)}{\sqrt x} \quad\text{and}\quad E^\pi(x;q,a) = \frac{\Delta^\pi(x;q,a)}{\sqrt x/\log x}.
\]
It is convenient at times to use a prime counting function itself as the main term, and such error terms are denoted by the symbol~$\Deltadot$: for example,
\[
\Deltadot^\psi(x;q,a) = \phi(q)\psi(x;q,a) - \psi(x) \quad\text{and}\quad \Deltadot^\pi(x;q,a) = \phi(q)\pi(x;q,a) - \pi(x).
\]
(Typically this modification results in the same explicit formula with the principal character removed.)
The corresponding normalized error terms are denoted by~$\Edot$: for example,
\[
\Edot^\psi(x;q,a) = \frac{\Deltadot^\psi(x;q,a)}{\sqrt x} \quad\text{and}\quad \Edot^\pi(x;q,a) = \frac{\Deltadot^\pi(x;q,a)}{\sqrt x/\log x}.
\]
We extend our convention regarding counting functions in arithmetic progressions: for example,
\[
\Delta^\psi(x;q,a,b) = \Delta^\psi(x;q,a) - \Delta^\psi(x;q,b) \quad\text{and}\quad E^\pi(x;q,a,b) = E^\pi(x;q,a) - E^\pi(x;q,b).
\]
Note that functions of the first type are almost redundant, since (for example) $\Delta^\psi(x;q,a,b) = \phi(q) \psi(x;q,a,b)$ exactly. (And recall that some authors use $\Delta$ to mean this difference function without the factor $\phi(q)$.) However, there will be situations where each notation is useful to us; furthermore, this new use of~$\Delta$ already follows from existing notational conventions.


It can also be convenient to define this notation for the function $\psi(x,\chi) = \sum_{n\le x} \Lambda(n) \chi(n)$, for any Dirichlet character~$\chi$ (see Section~\ref{Dirichlet section}), in the following way:
\[
\Delta^\psi(x,\chi) = \psi(x,\chi) - \begin{cases}
x, &\text{if } \chi=\chi_0, \\
0, &\text{if } \chi\ne\chi_0,
\end{cases}
\qquad
E^\psi(x,\chi) = \frac{\Delta^\psi(x,\chi)}{\sqrt x}.
\]

All of the functions in this section so far have been real-valued (except for the last paragraph where the functions are potentially complex-valued); in the context of primes in arithmetic progressions, it is often helpful to consider vector-valued functions. We use subscripts to indicate the modulus and residue classes---for example,
\[
\pi_{q;a_1,\dots,a_r}(x) = \bigl( \pi(x;q,a_1), \dots, \pi(x;q,a_r) \bigr)
\quad\text{and}\quad
\Edot^\psi_{q;a_1,\dots,a_r}(x) = \bigl( \Edot^\psi(x;q,a_1), \dots, \Edot^\psi(x;q,a_r) \bigr).
\]

\subsection{Weighted versions of prime counting functions}

It is common to vary these prime counting functions by attaching a weight to each term in the sum, changing for example $\sum_{n\le x} \Lambda(n)$ to $\sum_{n\le x} \Lambda(n) g(n)$. We use the following consistent notation for the most common of these variants.

As is standard, the subscript~$0$, as in the example $\psi_0(x) = \frac12\big( \psi(x-) + \psi(x+) \big)$, represents a modification of a function's value at a jump discontinuity to equal the average of the left- and right-hand limits.

The subscript $r$ represents weighting by a reciprocal factor (often resulting in a ``Mertens sum''); for example,
\[
\pi_r(x) = \sum_{p\le x} \frac 1p, \quad  \theta_r(x) = \sum_{p\le x} \frac{\log p}{p}, \quad\text{and}\quad \psi_r(x;q,a) = \sum_{\substack{n\le x \\ n\equiv a\mod q}} \frac{\Lambda(n)}n.
\]
If we wish to modify one of these Mertens sums at its jump discontinuities as above, we concatenate the two subscripts: for example, $\pi_{r0}(x) = \frac12\big( \pi_r(x-) + \pi_r(x+) \big)$. Indeed all of our previous notational variants can apply to these sums as well---for example,
\[
\Delta^{\pi_r}(x) = \pi_r(x) - (\log\log x + B) \quad\text{and}\quad E^{\pi_r}(x) = \sqrt x\log x \cdot \Delta^{\pi_r}(x)
\]
for the appropriate constant~$B$.

The subscript $e$ represents weighting by an exponentially decaying function of $x$ rather than cutting off abruptly at $x$; for example,
\[
\pi_e(x) = \sum_p e^{-p/x} \quad\text{and}\quad \psi_e(x;q,a) = \sum_{\substack{n\ge 1 \\ n\equiv a\mod q}} \Lambda(n)e^{-n/x}.
\]
In terms of their asymptotics, these exponentially weighted sums usually act like their abrupt-cutoff versions; for example, $\pi_e(x)$ has a similar size to $\pi(x)$. However, their oscillations are typically damped, often resulting in rather different properties when comparing two such functions to each other (such as the exponentially weighted version having a bias for one sign while the unweighted version exhibits oscillations of sign).

The subscript $l$ represents weighting by a certain exponential factor with a squared logarithm, scaled by a second parameter $r$: for example,
\[
\pi_l(x,r) = \sum_p e^{-\frac1r(\log \frac px)^2} \quad\text{and}\quad \psi_l(x,r;q,a) = \sum_{\substack{n\ge 1 \\ n\equiv a\mod q}} \Lambda(n)e^{-\frac1r(\log \frac nx)^2}.
\]
In asymptotic terms, this weighting is similar to restricting the range of summation to approximately $[e^{-\sqrt r} x, e^{\sqrt r} x]$; again, the oscillatory nature of the weighted sum can be rather different.

When the weight function is a Dirichlet character $\chi$ (see Section~\ref{Dirichlet section}), we follow the tradition of putting $\chi$ as an extra argument rather than a subscript; for example,
\[
\theta(x,\chi) = \sum_{p\le x} \chi(p)\log p.
\]

\subsection{Summatory functions}

Certain summatory functions of multiplicative functions have been analyzed using the techniques of comparative prime number theory. Two notable examples are the sums of the M\"obius and Liouville functions, which are denoted by
\[
M(x) = \sum_{n\le x} \mu(n) \quad\text{and}\quad L(x) = \sum_{n\le x} (-1)^{\Omega(n)},
\]
respectively. (The Liouville function is typically denoted by $\lambda(n) = (-1)^{\Omega(n)}$, but we avoid that notation herein to free the symbol $\lambda$ for other uses.) Two conjectures that motivated substantial work in comparative prime number theory are the ``Mertens conjecture'', the assertion that $|M(x)| < \sqrt x$, and the ``P\'olya problem'', the assertion that $L(x) \le 0$. (The latter assertion is often mistakenly named ``P\'olya's conjecture'', but P\'olya only posed and studied the problem rather than making a definitive conjecture and indeed probably found it unlikely to be true.) Both assertions have been disproved (in~\cite{MR783538} and~\cite{MR0104638}, respectively), although research continues into the distribution of these two functions. The weak Mertens conjecture, namely the assertion that $M(x) \ll \sqrt x$, is still unresolved, although it was shown in~\cite{MR0006202} to be incompatible with the pair of conjectures RH and LI (see Section~\ref{hypotheses section}).

The notational conventions from the previous sections are used for weighted versions of these summary functions as well; for example,
\[
M(x;q,a) = \sum_{\substack{n\le x \\ n\equiv a\mod q}} \mu(n) \quad\text{and}\quad L_r(x) = \sum_{n\le x} \frac{(-1)^{\Omega(n)}}{n};
\]
the conjecture that the latter is always nonnegative (often attributed to Tur\'an, though again he only studied the problem rather than asserting a conjecture) was also disproved in~\cite{MR0104638}.
We also define the notation $\Delta^M(x) = M(x)$ and $\Delta^L(x) = L(x)$ and $\Delta^{L_r}(x) = L_r(x)$; while unprofitable on their own, these definitions allow us to employ the notation for repeated averaging described in Section~\ref{errors section}, as well as the notation $E^M(x) = M(x)/\sqrt x$ and $E^L(x) = L(x)/\sqrt x$ and $E^{L_r}(x) = L_r(x)\sqrt x$.

We also introduce some standard notation for $k$-free numbers, which are numbers not divisible by the $k$th power of any prime, so that squarefree numbers are the case $k=2$. Let $Q_k(x)$ denote the number of $k$-free integers up to~$x$, and define $\Delta^{Q_k}(x)=Q_k(x)-x/\zeta(k)$. 
For integers $k\ge2$, the generalized M\"obius function $\mu_k(n)$ is defined to be $\mu_k(n)=(-1)^{\Omega(n)}$ if $n$ is $k$-free and $\mu_k(n)=0$ otherwise. Note that these functions interpolate between $\mu_2(n)=\mu(n)$ and $\lim_{k\to\infty} \mu_k(n)=(-1)^{\Omega(n)}$.
These quantities are related by the identity $Q_k(x)=\sum_{n \leq x}\mu_k^2(n)$.

Another summatory function studied using techniques that
overlap with those of comparative prime number theory is
\[
  D(x) = \sum_{n \le x} \tau(n) ,
  \text{ where } \tau(n) = \#\{ d \colon d \mid n \} = \sum_{d \mid n} 1 .
\]
It was first proven by Dirichlet (see~\cite{MR1576550} for a discussion of the history) that 
\[
  D(x) = x\log x + (2C_0-1)x + O(\sqrt x),
\]
where $C_0$ is Euler's constant. The study of the error
term $\Delta^D(x) = D(x) - x\log x - (2C_0-1)x$ is intertwined with comparative prime number theory, and one early result by Hardy \cite{MR1576550} demonstrates that the techniques of comparative prime theory are often applicable to the study of this error term.

\subsection{Counting sign changes}

We use the letter $W$ generally to denote the function that counts the number of sign changes of another function on an interval.
To be pedantic, if $h$ is a function from $(1,\infty)$ to $\R$, then we define
\begin{multline*}
W(h;T) = \max \big\{ n\ge0\colon \text{there exist } 1 < t_0 < t_1 < \cdots < t_n < T \\
\text{ with } h(t_{j-1}) h(t_j) < 0 \text{ for all } 1\le j\le n \big\}.
\end{multline*}
(One could quibble over whether taking the value $0$ counts as a sign change regardless of its neighboring values; the results in this subject tend not to require this loophole.)
We can demand large oscillations to go along with our sign changes by adding a function as an additional argument:
\begin{multline*}
W\big( h;T; S(t) \big) = \max \big\{ n\ge0\colon \text{there exist } 1 < t_0 < t_1 < \cdots < t_n < T \\
\text{ with } h(t_{j-1}) h(t_j) < 0 \text{ for all } 1\le j\le n \text{ and } |h(t_j)| > S(t_j) \text{ for all } 0\le j\le n \big\}.
\end{multline*}

Given functions $f$ and $g$ from $(1,\infty)$ to $\R$, we further define $W(f,g;T) = W(f-g;T)$ to be the counting function of sign changes of the difference $f(x)-g(x)$. Certain special cases of this notation deserve a shorthand notation: we define $W^\pi(T) = W(\pi,\li;T)$ and $W^\Pi(T) = W(\Pi,\li;T)$, and also $W^\theta(T) = W(\theta,x;T)$ and $W^\psi(T) = W(\psi,x;T)$ where $x$ denotes the identity function. (As before, we do not distinguish in our summaries between $\li(x)$ and $\Li(x)$ and $\sum_{2\le n\le x} 1/\log n$ in this context.) The bare notation $W(T)$ is a further shorthand for~$W^\pi(T)$.

In addition, given a positive integer $q$ and distinct reduced residues $a$ and $b\mod q$, we define $W^\psi_{q;a,b}(T) = W \big( \psi(x;q,a),\psi(x;q,b);T \big)$, and similarly with $\psi$ replaced by $\theta$, $\Pi$, or $\pi$; we further shorten $W^\pi_{q;a,b}(T)$ to $W_{q;a,b}(T)$. We may add a function as an additional argument as above to indicate large oscillations, as in $W_{q;a,b}(T;S(t))$ for example; similarly, we may replace single residue classes with sets of residue classes, as in $W_{q;\cN,\cR}(T)$.

\subsection{Densities}

The natural density of a set $\cS$ of positive real numbers is
\[
\mathfrak d(\cS) = \lim_{x\to\infty} \frac{\mathop{\rm meas}\big( \{ 0<t\le x\colon t\in \cS \} \big) }x = \lim_{x\to\infty} \frac1x \int\limits_{\substack{0<t<x \\ t\in \cS}} dt,
\]
where ``meas'' denotes Lebesgue measure on~$\R$.
On the other hand, the logarithmic density of a set $\cS \subset (1,\infty)$ is
\[
\delta(\cS) = \lim_{x\to\infty} \frac1{\log x} \int\limits_{\substack{1<t<x \\ t\in \cS}} \frac{dt}t.
\]
An easy change of variables shows that the logarithmic density of $\cS$ equals the natural density of the set $\log \cS = \{ \log t\colon t\in \cS \}$. Moreover, a partial summation argument shows that if the natural density $\mathfrak d(\cS)$ exists, then the logarithmic density $\delta(\cS)$ also exists and has the same value. However, there are sets whose natural density does not exist but whose logarithmic density does exist; for example, the union (over $k\in\N$) of the intervals $[10^{2k-1},10^{2k})$ has logarithmic density equal to~$\frac12$ but does not have a natural density.

We will use many variants of this logarithmic density notation. If $f_1,\ldots,f_r$ are functions from $(1,\infty)$ to~$\R$, then we define the shorthand notation
\[
\delta(f_1,f_2,\ldots,f_r) = \delta\big( \{ x>1\colon f_1(x) > f_2(x) > \cdots > f_r(x) \} \big).
\]
For example, $\delta(\li,\pi)$ is the logarithmic density of the set of real numbers $x>1$ for which $\li(x) > \pi(x)$. Certain special cases of this notation can be even further abbreviated. For example, let $q$ be a positive integer, and let $a_1,\ldots,a_r$ be distinct reduced residues\mod q. Then we define
\[
\delta_{q;a_1,\ldots,a_r} = \delta\big( \pi(x;q,a_1), \ldots, \pi(x;q,a_r) \big) = \delta\big( \{ x>1\colon \pi(x;q,a_1) > \cdots > \pi(x;q,a_r) \} \big).
\]
We also define
\[
\delta_{q;\cN,\cR} = \delta\big( \pi(x;q,\cN), \pi(x;q,\cR) \big) = \delta\big( \{ x>1\colon \pi(x;q,\cN) > \pi(x;q,\cR) \} \big)
\]
and similarly for $\delta_{q;\cR,\cN}$ (these definitions are sensible when $q$ has primitive roots).

Finally, we define the upper and lower logarithmic densities of $\cS$ (which always exist) as
\[
\overline\delta(\cS) = \limsup_{x\to\infty} \frac1{\log x} \int\limits_{\substack{1<t<x \\ t\in \cS}} \frac{dt}t, \quad
\underline\delta(\cS) = \liminf_{x\to\infty} \frac1{\log x} \int\limits_{\substack{1<t<x \\ t\in \cS}} \frac{dt}t,
\]
so that $\delta(\cS)$ exists if and only if $\overline\delta(\cS) = \underline\delta(\cS)$. This notation propagates through our shorthand notations as well; for instance, $\underline\delta_{q;\cN,\cR} = \underline\delta\big( \{ x>1\colon \pi(x;q,\cN) > \pi(x;q,\cR) \} \big)$.

\subsection{Limiting distributions and density functions}

Given a function $h\colon [0,\infty)\to\R$, the limiting (or asymptotic) cumulative distribution function of $h$ is the nondecreasing function
\[
\lim_{T\to\infty} \frac {\mathop{\rm meas}\{ t\in[0,T]\colon h(t) \le a \}}T = \lim_{T\to\infty} \biggl( \frac1T \int\limits_{\substack{0\le t\le T \\ h(t) \le a}} \,dt \biggr)
\]
if the limit exists (except at jump discontinuities, of which there are only a countable number). More common in comparative prime number theory is the limiting logarithmic cumulative distribution function, with the analogous definition
\[
\kappa^h(\alpha) = \lim_{U\to\infty} \frac1{\log U} \biggl( \int\limits_{\substack{1\le u\le U \\ h(u) \le \alpha}} \,\frac{du}u \biggr),
\]
which equivalently is the cumulative distribution function of $h(e^t)$. There is a corresponding limiting logarithmic density $\mu^h$, which is the measure satisfying
\[
\mu^h\bigl( (\alpha,\beta] \bigr) = \int_\alpha^\beta \kappa^h(x) \,dx
\]
for any real numbers $\alpha<\beta$. It has the property that for any bounded continuous function $f(x)$,
\[
\lim_{U\to\infty} \frac1{\log U} \biggl( \int_1^U f\bigl( h(u) \bigr) \,\frac{du}u \biggr) = \int_{\R} f(x)\,d\mu^h(x),
\]
and the continuity assumption can be omitted if $\mu^h$ is absolutely continuous with respect to Lebesgue measure. These logarithmic densities are probability measures and thus can be viewed as the densities of random variables. Vector-valued functions have analogous logarithmic cumulative distribution functions and logarithmic densities on~$\R^r$.

\section{Notation related to complex analysis} \label{other notation section}

As is usual in analytic number theory, we often use $s=\sigma+it$ to denote a complex variable and its real and imaginary parts; its argument will be denoted by $\arg(s)$, so that $s=|s|e^{i\arg s}$. If $\rho$ is a nontrivial zero of a Dirichlet (or other) $L$-function, including the Riemann zeta-function, we write $\rho=\beta+i\gamma$ to refer to its real and imaginary parts.

\subsection{Dirichlet characters and Dirichlet $L$-functions} \label{Dirichlet section}

As usual, a Dirichlet character with modulus $q$ is a completely multiplicative function on $\Z$ with period $q$ whose support is the set of integers coprime to~$q$. We call characters real, complex, quadratic, (im)primitive, and induced with their standard meanings; the conductor of a character $\chi$ is the modulus of the primitive character $\chi^*$ that induces it.

We use $\chi_0$ to denote the principal character (the modulus being understood from context). When $D\ne1$ is a fundamental discriminant, we let $\chi_D$ denote the associated quadratic character, which is a primitive character of conductor~$|D|$ that is even if $D>0$ and odd if $D<0$. When $q$ is prime, we use the shorthand $\chi_{\pm q}$ to mean $\chi_q$ if $q\equiv1\mod4$ and $\chi_{-q}$ if $q\equiv3\mod4$. On the other hand, by $\chi_1$ we mean a hypothetical quadratic character with an exceptional zero $\beta_1$.

Every Dirichlet character gives rise to a Dirichlet $L$-function $L(s,\chi) = \sum_{n=1}^\infty \chi(n)n^{-s}$. Like the Riemann zeta-function (which is the special case $q=1$ and $\chi=\chi_0$), Dirichlet $L$-functions have infinitely many nontrivial zeros $\rho=\beta+i\gamma$ in the critical strip $0<\beta<1$. These zeros are counted by the function
\[
N(T,\chi) = \#\{ \rho\colon L(\rho,\chi) = 0,\, 0<\beta<1,\, |\gamma|\le T\}.
\]
Note the slight dissonance with the traditional notation
\[
N(T) = \#\{ \rho\colon \zeta(\rho) = 0,\, 0<\beta<1,\, 0\le \gamma\le T\}
\]
which counts only nontrivial zeros of $\zeta(s)$ in the upper half-plane: this suffices for $\zeta(s)$ due to the Schwarz reflection principle, but Dirichlet $L$-functions do not all possess that symmetry.

Sums over zeros of Dirichlet $L$-functions (of the type that arise in explicit formulas, for example) often do not converge absolutely, and therefore we adopt the standing convention that sums over nontrivial zeros are limits of their symmetric truncations:
\[
\sum_\rho f(\rho) = \lim_{T\to\infty} \sum_{\substack{L(\rho,\chi)=0 \\ 0<\beta<1 \\ |\gamma|\le T}} f(\rho).
\]

\subsection{Landau's theorem}

For a real-valued function $A(x)$, define
\[
g(s)=\int_1^\infty \frac {A(x)}{x^s} \,dx 
\]
Typically there will be a real number $\sigma_0$ such that this integral converges when $\sigma>\sigma_0$ and diverges when $\sigma<\sigma_0$. 
Landau proved that if $A(x)$ is eventually positive or eventually negative, then $g(s)$ has a singularity at $s=\sigma_0$ (that is, $g(s)$ must have a rightmost singularity on the real axis).

The contrapositive of this theorem is a useful tool in comparative prime number theory:
Suppose that $g(s)$ has no singularities on the subray $\{ \sigma\in\R\colon \sigma>\sigma_1\}$ of the real axis (that is, $g(s)$ is analytic on a neighborhood of that ray), but $g(s)$ is not analytic in the half-plane $\{s\in\C\colon \sigma>\sigma_1\}$. Then $A(x)$ has arbitrarily large sign changes.

\subsection{Explicit formulas}\label{formulas section}

As mentioned earlier, one of the defining characteristics of comparative prime number theory is the presence of an ``explicit formula''. There is no precise definition of that term, but typically an explicit formula contains a sum over the (nontrivial) zeros of some $L$-function. The prototypical example is the explicit formula
\[
\psi_0(x) = x - \sum_{\rho} \frac{x^\rho}\rho - \log2\pi - \frac12\log\bigg( 1 - \frac1{x^2} \bigg)
\]
for the Chebyshev function $\psi(x)$ modified at its jump discontinuities; the fact that this is an exact equality for all $x>1$ is one of the most beautiful statements in analytic number theory.

Explicit formulas for prime-counting functions yield explicit formulas for their normalized error terms; for example, assuming the generalized Riemann hypothesis,
\[
E^\theta(x;q,a,b) = c_q(b)-c_q(a) - \sum_{\chi\mod q} \big( \overline\chi(a) - \overline\chi(b) \big) \sum_{\substack{\gamma\in\R \\ L(1/2+i\gamma,\chi) = 0}} \frac{x^{i\gamma}}{\frac12+i\gamma} + O_q(x^{-1/6}).
\]
This formula is helpful for studying when $E^\theta(x;q,a,b) > 0$, or equivalently when $\theta(x;q,a) > \theta(x;q,b)$. Note that each summand in the inner sum oscillates around a circle of fixed radius (one that decreases as~$\gamma$ increases); while this inner sum is not literally bounded, it is bounded on average over~$x$ and possesses a limiting logarithmic distribution. Therefore $E^\theta(x;q,a,b)$ has a limiting logarithmic distribution with mean $c_q(b)-c_q(a)$, the sign of which depends on whether~$a$ and~$b$ are quadratic residues or nonresidues modulo~$q$.

\subsection{The power-sum method}

A great deal of early progress in comparative prime number theory, particularly the unconditional results, relied on the study of linear combinations of powers of complex numbers, namely sums of the shape
\[
s_{v}=\sum_{j=1}^{n} b_{j} z_{j}^{v}.
\]
Lower bounds for such sums were systematically developed by Tur\'an and S\'os. While there are many variants of these lower bounds that have been obtained, they can be grouped into two main categories.

The ``first main theorem'' is a type of result that applies when the $z_j$ are large. For example,
suppose that $z_1,\ldots,z_n$ are distinct complex numbers with $|z_n|\ge1$ for all~$n$. For any nonnegative integer $m$, there exists an integer $m+1\le v\le m+n$ such that
\begin{equation*}
|s_{v}| \ge \bigg(\frac{n}{A(m+n)}\bigg)^{n} |s_0|,
\end{equation*}
where $A$ is an absolute constant. 

The ``second main theorem'' is a type of result that applies when the $z_j$ are small. For example,
suppose that $z_1,\ldots,z_n$ are distinct complex numbers with $1\ge |z_1|\ge\cdots\ge|z_n|$. For any nonnegative integer $m$, there exists an integer $m+1\le v\le m+n$ such that
\begin{equation*}
|s_{v}| \ge \bigg(\frac{n}{B(m+n)}\bigg)^{n} \min _{1 \leq j \leq n} \bigg| \sum_{n=1}^{j} b_{n} \bigg|,
\end{equation*}
where $B$ is an absolute constant.

Instead of restricting the candidate exponents~$v$ to an interval of exactly $n$ consecutive integers, we may allow candidates from a longer range of exponents.
For example, in the ``second main theorem'' (so that $1\ge |z_1|\ge\cdots\ge|z_n|$), let $m\ge N\ge n$; then there exists an integer $m+1\le v\le m+N$ such that
\begin{equation*}
|s_{v}| \ge \bigg(\frac{N}{Bm}\bigg)^{N} \min_{1 \leq j \leq n} \bigg| \sum_{n=1}^{j} b_{n} \bigg|.
\end{equation*}
For the ``second main theorem'', one can also obtain better conclusions by adding an ``argument restriction'', that is, the assumption that each $|\arg z_j|\ge\ep$ for some fixed $\ep>0$. Stronger results can also be obtained by assuming that each $b_j$ is a nonnegative real number, and strengthened further by restricting to the special case $b_1=\cdots=b_n=1$.

Note that these results show that some $s_v$ is large in modulus but gives no information about its argument. Tur\'an (somewhat unhelpfully) calls these results ``two-sided'' theorems. There exist analogous results where the lower bound applies not just to $|s_v|$ but to $\Re s_v$ or $-\Re s_v$; Tur\'an calls such results ``one-sided'' theorems.

\subsection{$k$-functions}  \label{k-functions section}

A great deal of the work of Kaczorowski involves certain functions called $k$-functions, which are superficially similar to sums that appear in explicit formulas for $\psi(x,\chi)$. For $\Im z>0$, define
\begin{align*}
k(z, \chi) = \sum_{\gamma>0} e^{\rho z} \quad \text{and} \quad K(z, \chi) = \sum_{\gamma>0} \frac{e^{\rho z}}{\rho},
\end{align*}
where the sums are over zeros of $L(s,\chi)$ in the upper half-plane.

These functions can be regarded as having their domain equal to~$\cM$, the Riemann surface for $\log z$; every point on the surface can be uniquely written as $re^{ia}$ where $r>0$ and $a\in\R$. Let $z^c$ denote the natural extension of complex conjugation to~$\cM$, namely $(re^{ia})^c = re^{-ia}$; also let $z^*$ denote an extension of multiplication by $-1$ to~$\cM$, namely $(re^{ia})^* = re^{i(a - \pi)}$.

Certain functions appear frequently in connection to $k$-functions: define
\[
D(z, \chi) = - \sum_{\substack{\beta>0 \\ L(\beta, \chi)=0}} e^{\beta z} + \frac1{e^{2z}-1} \begin{cases}
e^{3z} + e^{2z}-1, &\text{if } \chi=\chi_0, \\
e^z, &\text{if } \chi\ne\chi_0 \text{ and } \chi(-1)=1, \\
e^{2z}, &\text{if } \chi(-1) = -1.
\end{cases}
\]
Further define
\begin{align*}
F(x,\chi) = \lim_{y \to 0^+} \bigg( K(x+iy, \chi) + \overline{K(x+iy, \overline{\chi})} \bigg)
\end{align*}
and
\[
R_1(x) = \frac12\log(1-e^{-2x}), \quad R_{-1}(x) = \frac12 \log \frac{e^x-1}{e^x+1}.
\]
Certain constants also appear frequently: define
\[
B(\chi) = \sum_{\substack{\beta>0 \\ L(\beta, \chi)=0}}  \frac{1}{\beta} - \frac{C_0}{2} -\frac{1}{2} \log \frac{\pi}{q} + F(0, \chi) - \begin{cases}
1, &\text{if } \chi=\chi_0, \\
0, &\text{if } \chi\ne\chi_0 \text{ and } \chi(-1)=1, \\
\log2, &\text{if } \chi(-1) = -1
\end{cases}
\]
(note that $B(\chi)$ is not the same as a constant of the same name related to the Hadamard product expansion of $L(s,\chi)$) and $C(\chi) = B(\chi) + C_0 + \log\frac{2\pi}q$.

\subsection{Hypotheses on zeros}\label{hypotheses section}

It is extremely difficult to obtain unconditional results in comparative prime number theory, particularly where limiting logarithmic distributions and densities are concerned. Certain assumptions on the zeros of Dirichlet $L$-functions therefore arise repeatedly in this subject. The most famous of these is the generalized Riemann hypothesis (GRH), sometimes called the Riemann--Piltz conjecture, which asserts that all nontrivial zeros of all Dirichlet $L$-functions have real part equal to~$\frac12$. We use $\sigma_0$-GRH to denote the weaker (but still currently inaccessible) assertion that $L(\sigma+it,\chi)$ does not vanish when $\sigma>\sigma_0$, so that $1$-GRH is trivial and $\frac12$-GRH is the same as the full GRH.

Given a nonempty set $X$ of Dirichlet $L$-functions (or, abusing notation slightly, Dirichlet characters), we let $\Theta(X)$ denote the supremum of the real parts of their zeros, that is, the smallest real number such that $\Theta(X)$-GRH holds. We use the abbreviation $\Theta(q)$ when~$X$ is the set of all Dirichlet characters modulo~$q$, as well as $\Theta(\chi)$ when~$X$ consists of the single Dirichlet character~$\chi$. The assertion that some Dirichlet $L$-function in~$X$ has a zero with real part exactly equal to $\Theta(X)$ is abbreviated SA for ``supremum attained'' (and sometimes referred to as ``Ingham's condition''). We may write SA$(X)$ to emphasize that we are considering a specific set of Dirichlet $L$-functions, but the set is often inferred from context (this remark applies similarly to the remainder of the notation in this section). We note that GRH implies SA but that $\Theta(X)=1$ is inconsistent with~SA.

Regarding the vertical distributions of the zeros, we use HC to denote the ``Haselgrove condition'' that no Dirichlet character (in the set under discussion) vanishes on the segment $0<\sigma<1$ of the real axis. Such a real zero would create a non-oscillatory term in relevant explicit formulas, one that could result in an unexpected source of bias. By continuity, HC implies that there exists a positive constant~$E_k$ such that these $L(s,\chi)$ are nonzero on the rectangle $\{ 0<\sigma<1,\, |t| \le E_k \}$; we write HC$(E_k)$ if we need to refer to this parameter.

The notation GRH$(H)$ (sometimes called the ``finite Riemann--Piltz'' conjecture) denotes the generalized Riemann hypothesis ``up to height~$H$'', namely the statement that if $\rho$ is a nontrivial zero of $L(s,\chi)$ with $|\gamma|\le H$, then $\beta=\frac12$. 
Note that HC$(E_k)$ implies GRH$(H)$ if $E_k \ge H$; on the other hand, GRH$(H)$ gives no constraint at all upon zeros on the critical line. We therefore use the notation GRH$(H,E_k)$ to denote the combination of GRH$(H)$ and HC$(E_k)$, the latter of which constrains only the zeros on the critical line when $E_k < H$.
Note also that GRH$(0)$ is almost the same as HC, except that GRH$(0)$ allows for the possibility of a zero at $s=\frac12$.

The arithmetic nature of the imaginary parts (ordinates) of zeros of $L(s,\chi)$ is also significant in comparative prime number theory. We write LI (sometimes called GSH for the ``grand simplicity hypothesis'') to denote the ``linear independence'' assertion that the multiset of nonnegative ordinates of zeros of the relevant Dirichlet $L$-functions is linearly independent over the rational numbers. In particular, LI implies that all zeros are simple and that $L(\frac12,\chi)\ne0$. We use LI$(\sigma)$ to denote the corresponding linear independence conjecture restricted to the zeros with real parts greater than or equal to~$\sigma$.

For the Riemann zeta-function, the Riemann hypothesis (RH) is the assertion that all nontrivial zeros of $\zeta(s)$ have real part equal to~$\frac12$. Almost all of the other notation above would be used in the same form when referring to $\zeta(s)$, although $\Theta(\{\zeta(s)\})$ is abbreviated simply to~$\Theta$. These same abbreviations are also used for analogous hypotheses on zeros of other $L$-functions that should be clear from context.

\section{Types of questions} \label{questions section}

Given two functions $f,g\colon (1,\infty)\to\R$ that are asymptotic to each other, such as $\pi(x)$ and $\li(x)$ or $\pi(x;4,1)$ and $\pi(x;4,3)$, the questions that comparative prime number theory tends to ask about the pair of functions are:
\begin{enumerate}
\item Are there arbitrarily large values of $x$ for which $f(x) > g(x)$, and arbitrarily large values of $x$ for which $g(x) < f(x)$? In other words, does the difference $f(x)-g(x)$ change signs infinitely often? (These are not quite mathematically identical because of the possibility of plentiful or carefully arranged ties $f(x)=g(x)$, so implicit in this question is asking whether such ties are rare.) The other alternative is that one of the functions exceeds the other for all sufficiently large~$x$.
\item How large and positive can the difference $f(x)-g(x)$ get? How large and negative can it get?
\item More generally, what is the distribution of values of $f(x)-g(x)$? Is it possible that some suitably normalized version of this difference, such as $(f(x)-g(x))/\sqrt x$, actually has a limiting distribution or a limiting logarithmic distribution?
\item How often does the difference $f(x)-g(x)$ change sign? How many sign changes are there in $(1,X)$ as a function of $X$? How close can we take $Y=Y(X)$ to $X$ to ensure that there is always a sign change in $[X,Y]$?
\item What is the natural density of the set of real numbers $x>1$ for which $f(x)>g(x)$? What is its logarithmic density~$\delta(f,g)$? (Typically we believe that the natural densities of such sets do not exist in prime number races, but that their logarithmic densities do exist.)
\item Given a family of races, such as $\pi(x;q,\cN)$ versus $\pi(x;q,\cR)$: how do answers to the above questions, such as $\delta_{q;\cN,\cR}$, depend upon the member of the family ($q$ in this case)? Do the distributions of the members of the family tend to some limit, such as a normal distribution?
\end{enumerate}
Some of the above questions have analogues for several functions $f_1,\ldots,f_r\colon(1,\infty)\to\R$ considered together:
\begin{enumerate}\setcounter{enumi}{6}
\item Are there arbitrarily large values of $x$ for which $f_1(x) > \cdots > f_r(x)$? Does this remain true no matter how we permute the~$f_j$?
\item More generally, what is the distribution of values of the vector $\big( f_1(x), \ldots, f_r(x) \big) \in \R^r$? Is it possible that some suitably normalized version of this difference actually has a limiting distribution or a limiting logarithmic distribution?
\item What is the natural density of the set of real numbers $x>1$ for which $f_1(x) > \cdots > f_r(x)$? What is its logarithmic density $\delta(f_1,\ldots,f_r)$? (As before, we believe that the natural densities of such sets do not exist in prime number races, but that their logarithmic densities do exist.)
\item Given a family of such $r$-way races, how do answers to the above questions depend upon the member of the family? Do the distributions of the members of the family tend to some limit, such as a multivariate normal distribution?
\end{enumerate}

The articles~\cite{MR0146156} and~\cite{MR0162771} by Knapowski and Tur\'an present organized schema for problems in comparative prime number theory, as do surveys of these topics by Kaczorowski~\cite{MR1645811} and by Ford and Konyagin~\cite{MR1985941}, although several of the questions listed above had not yet been investigated sufficiently deeply to make some of their lists.

\section*{Acknowledgments}

We gratefully thank Devang Agarwal,
Alexandre Bailleul,
Michael Coons,
Alia Hamieh,
Elchin Hasanalizade,
Daniel R.~Johnston,
Farid Jokar,
Florent Jouve,
Shin-ya Koyama,
\LaTeX\ Stack Exchange user ``moewe'',
Michael J.~Mossinghoff,
Nathan Ng, and
Alan Xiang
for their contributions to this manuscript. We also thank the anonymous referees for their thorough readings and detailed suggestions for corrections and improvements.
Many authors' research was supported by the Natural Science and Engineering Research Council of Canada.


\section*{Chronological bibliography}

The annotated bibliography begins here, with all of the sources cited and summarized listed in chronological order; items in this chronological list are labeled by their number alone, such as~[123]. Following the annotated bibliography is a second list, in alphabetical order by author, of the same set of sources but without annotations; items in this alphabetical list have been given labels that are numbers following the letter ``A'' (for ``alphabetical''), such as~[A45], to distinguish them from the labels in the main list. Each entry in the second bibliography links to its corresponding entry and annotation in the first bibliography.

Our goal has been to describe the results using a single system of notation, both to avoid the need to define notation in individual annotations and to propose a unified notation for current and future practitioners of comparative prime number theory. Any notation in a summary that is not defined there can be found or deduced from the detailed material in Sections~2--3.


\bigskip
\label{first bib}
   
\normalsize
\printbibliography[heading=none]

\section*{Alphabetic bibliography} \label{second bib}
\setlength{\bibitemsep}{-.2\baselineskip}
\togglefalse{showannotation}
{
\localrefcontext[sorting=nyt, labelprefix=A]
\printbibliography[heading=none]
}

\end{document}